\documentclass[a4paper,11pt,twoside]{article}

\usepackage{amssymb,amsmath,amsfonts, amsthm,amscd}

\usepackage{verbatim}

\input xy
\xyoption{all}

\usepackage{style}

\newcommand{\Lie}{\mathfrak}

\title{Twisted cohomology for hyperbolic three
manifolds.}
\author{Pere Menal-Ferrer \and Joan Porti 
\thanks{Both authors partially supported by the 
Spanish Micinn through grant MTM2009-0759 and by 
the Catalan AGAUR through grant SGR2009-1207. 
The second author received the prize “ICREA Acad\`emia” 
for excellence in research, funded by the Generalitat de Catalunya.}
}

\begin{document}
\maketitle

\begin{abstract}
For a complete hyperbolic three manifold $M$, we consider the representations
of $\pi_1(M)$ obtained by composing a lift of the holonomy with complex
finite dimensional representations of $\operatorname{SL}(2,\mathbf C)$.
We prove a vanishing result for the cohomology of $M$ with coefficients 
twisted by these representations, using techniques of Matsushima-Murakami.
We give some applications to local rigidity.
\end{abstract}


Let $M$ be an orientable complete hyperbolic three manifold.
The holonomy representation of the complete hyperbolic structure 
$$\operatorname{Hol}\colon \pi_1(M)\to\operatorname{Isom}^+\,\mathbf{H}^3 \cong 
\operatorname{PSL}(2,\mathbf{C}),$$
can be lifted to a representation 
$\widetilde{\operatorname{Hol}}\colon \pi_1(M) \rightarrow \operatorname{SL}(2,\mathbf{C})$ 
(see for instance \cite{culler}), and there is a one-to-one correspondence 
between these lifts and spin structures on $M$. Composing one of these lifts
with a finite representation $V$ of $\operatorname{SL}(2,\mathbf{C})$,
we obtain a representation $\rho \colon \pi_1 (M) \rightarrow \operatorname{SL}(V)$.
Then we can consider the associated flat vector bundle $E_\rho$.

We will consider only \emph{complex} and \emph{finite dimensional} representations of
$\operatorname{SL}(2,\mathbf{C})$. It is well known that for
every positive integer $n$ there exists only one complex irreducible
representation $V_n$ of $\operatorname{SL}(2,\mathbf{C})$ of dimension $n$.
Moreover, $V_n$ is $(n-1)$-th symmetric power of the standard representation
$V_2 = \mathbf{C}^2$. Let\[
	\rho_n \colon \pi_1(M) \rightarrow \operatorname{SL}(n, \mathbf{C}).
\] 
denote the representation 
$\rho$ defined above for $V_n$.

A hyperbolic $3$-manifold $M$ is said to be \emph{topologically finite} 
if it is the interior of a compact manifold $\overline M$.
This is equivalent to say that $\pi_1(M)$ is finitely generated, by 
the proof of Marden's conjecture 
\cite{Agol, CalGab}.

Along the paper we shall assume that $M$ is \emph{nonelementary}, which means, in the context of three manifolds,
that its holonomy is an irreducible representation in $\operatorname{PSL}(2,\mathbf{C})$, namely that there is 
no proper invariant subspace of $\mathbf C^2$. Elementary manifolds have a simple geometry and topology 
(cf.\ Lemma~\ref{lemma:elementary}) and the following 
results still hold and have a straightforward proof.

\begin{theorem}
	\label{thm:inclusion}
	Let $M$ be a complete, nonelementary, hyperbolic $3$-manifold that is topologically finite,
	and $n \geq 2$.
	Then the inclusion $\partial\overline{M}\subset \overline M$ induces
	an injection,	
	\[
	H^1(M; E_{\rho_n})\hookrightarrow H^1(\partial\overline M; E_{\rho_n}),
	\]
	with $\dim H^1(M; E_{\rho_n}) =\frac12 \dim H^1(\partial\overline M; E_{\rho_n})$, and
	an isomorphism 
	\[
	H^2(M; E_{\rho_n})\cong H^2(\partial\overline M; E_{\rho_n}).
	\]
\end{theorem}

If $M$ is a complete hyperbolic $3$-manifold of finite volume with a single cusp,
then $\partial \overline M$ is a torus. An analysis of the cohomology groups
$H^*( \partial \overline M ; E_\rho)$ shows that all these groups
vanish for the representations $\rho_{2k}$, with $k>0$ (see Section~\ref{sec:cohom_ends}).
Hence, using Theorem~\ref{thm:inclusion} we get the following result.

\begin{theorem} 
	\label{thm:H1cusp}
	Let $M$ be a complete hyperbolic $3$-manifold of \emph{finite volume}
	with a \emph{single cusp}. Then for $k \geq 1$ we have 
	\[
	H^*(M; E_{\rho_{2k}}) = 0.
	\]
\end{theorem}

Notice that this theorem applies to hyperbolic knot exteriors in $S^3$. For instance, it allows to compute
Reidemeister torsions for hyperbolic knot exteriors.

Theorem \ref{thm:inclusion} has applications to infinitesimal rigidity. 
The space of infinitesimal deformations of $\rho_n$
is isomorphic to $H^1(M; E_{\operatorname{Ad}\circ \rho_n} )$, where 
\[
\operatorname{Ad} \colon \operatorname{SL}(n,\mathbf{C}) \rightarrow \operatorname{Aut}( \Lie{sl}(n,\mathbf C) )
\]
is the adjoint representation. 

The following theorem is an infinitesimal rigidity result for
$\rho_n$ in $\operatorname{SL}(n,\mathbf C)$ relative to the boundary. 
Its proof uses the decomposition of representation $\Lie{sl}(n,\mathbf C)$ 
into irreducible factors, and will be given in Section \ref{section:infinitesimalrigidity}.

\begin{theorem}
	\label{thm:l2rigidity}
	Let $M$ be a complete, hyperbolic, nonelementary and orientable $3$-manifold that is topologically finite. 
	If $\partial \overline M$ is the union of $k$ tori and $l$ surfaces of genus $g_1,\ldots,g_l\geq 2$,
	and $n\geq 2$, then 
	\[
	\dim_\mathbf{C} H^1(M; E_{\operatorname{Ad}\circ \rho_n} ) = k(n-1) + \sum (g_i-1)(n^2-1).
	\]
	In particular, if $M$ is closed then $H^1(M; E_{\operatorname{Ad}\circ \rho_n}) = 0.$
	In addition, all nontrivial elements in $H^1(M; E_{\operatorname{Ad}\circ \rho_n} )$
	are nontrivial in $H^1(\partial\overline M; E_{\operatorname{Ad}\circ \rho_n} )$ and have no $L^2$ representative.
\end{theorem}

When $n = 2$, this is Weil's infinitesimal rigidity in the compact case, and 
Garland's $L^2$-infinitesimal rigidity in the noncompact case.
This has been generalized to cone three manifolds by Hodgson-Kerckhoff \cite{HKJdG},
Weiss \cite{Weiss} and Bromberg \cite{Bromberg}.

Let $X(M,\operatorname{SL}(n, \mathbf C))$ be the variety of  
characters of $\pi_1(M)$ in $\operatorname{SL}(n,\mathbf C)$. The character of $\rho_n$
is denoted by $\chi_{\rho_n}$.
From the previous theorem and standard results on the variety of characters, we deduce:

\begin{theorem}
	\label{thm:smoothness}
	Let $M$ be a topologically finite, hyperbolic, nonelementary and  orientable $3$-manifold as in Theorem~\ref{thm:l2rigidity}.
	If $n\geq 2$, 
	then the character $\chi_{\rho_n}$ is a smooth point of 
	$X(M,\operatorname{SL}(n,\mathbf C))$ with tangent space $H^1(M; E_{\operatorname{Ad}\circ \rho_n} ) $.
\end{theorem}

For $n=2$, this is Theorem~8.44 of Kapovich \cite{Kapovich}.

This paper is organized as follows. 
In Section~\ref{section:finitedimreps} we recall some results about finite dimensional complex representations of $\operatorname{SL}(2,\mathbf C)$. 
Section~\ref{sec:vanishing} is devoted to Raghunathan's vanishing theorem, from which Theorem~\ref{thm:inclusion} will follow. Theorem~\ref{thm:H1cusp}
is proved in Section~\ref{sec:cohom_ends_lifts}, where we compute the cohomology of the ends and discuss some
properties of lifts of representations.
Section~\ref{section:infinitesimalrigidity} deals with applications to infinitesimal and
local rigidity, in particular we prove Theorems~\ref{thm:l2rigidity}
and \ref{thm:smoothness}. 

Appendix~\ref{principalbundles} reviews some results about principal bundles that are required in Section \ref{sec:vanishing}.

\medskip
\section{Finite dimensional complex representations of $\operatorname{SL}(2,\mathbf C)$ }
\label{section:finitedimreps}

Irreducible complex finite dimensional representation of $\operatorname{SL}(2,\mathbf C)$
are well known to be the symmetric powers of the standard representation $\mathbf{C}^2$.
Therefore, there is exactly one irreducible representation in each dimension. 
Let $V_n$ denote the irreducible complex $n$-th
dimensional representation of $\mathbf{C}^2$. We have $V_{n} = \operatorname{Sym}^{n-1} V_2$,
with the convention that $\operatorname{Sym}^0$ is the base field.

The decomposition into irreducible factors of the tensor product of two given complex
irreducible representation is given by the Clebsch-Gordan formula (cf.\ \cite[\S 11.2]{FultonHarris}).

\begin{theorem}[Clebsch-Gordan theorem]
	\label{thm:CG}
	For non-negative integer numbers $n, k$ we have
	\[
	V_{n} \otimes V_{n+k} = \bigoplus_{i = 0}^{n-1} V_{2(n-i) + k - 1}.
	\]
\end{theorem}

\begin{lemma} 
	\label{lemma:pairing}
	Let $V$ a finite dimensional complex representation of $\operatorname{SL}(2,\mathbf C)$. 
	Then there exists a nondegenerate $\mathbf C$-bilinear invariant pairing 
	\[
	\phi\colon V \times V \to\mathbf C.
	\] 
	Moreover, if $V$ is irreducible this pairing is unique up to multiplication by nonzero 
	scalars.
\end{lemma}

\begin{proof}
	From the classification of the irreducible representations 
	of $\operatorname{SL}(2,\mathbf C)$, we deduce that
	$V^*$ is isomorphic to $V$. 
	Thus we get an invariant bilinear pairing
	by composing the isomorphism $V\times V\cong V^*\times V$ 
	with the natural pairing between $V^*$ and $V$. 
	If $V$ is irreducible, $V=V_n$, then  the Clebsch-Gordan formula implies
	that $(V_n\otimes V_n)^*\cong V_n\otimes V_n$ has only one irreducible factor of
	dimension $1$, so the bilinear pairing is unique in this case.
\end{proof}

From this lemma we get (cf. \cite[Sec. 2.2]{goldman}): 
\begin{corollary}
	\label{corollary:PD}
	Poincar\'e duality with coefficients in $E_{\rho}$ holds true. 
\end{corollary}

Let $\operatorname{Ad} \colon \operatorname{SL}(n,\mathbf{C}) \to \operatorname{Aut}( \Lie{sl}(n,\mathbf C) )$ 
denote
the adjoint representation of $\operatorname{SL}(n,\mathbf{C})$. Composing it with the representation $V_n$
we get a representation $\operatorname{SL}(2,\mathbf{C}) \to \operatorname{Aut}( \Lie{sl}(n,\mathbf C) )$, which
makes $\Lie{sl}(n,\mathbf C)$ a $\operatorname{SL}(2,\mathbf{C})$-module. Next we want to decompose this module
into irreducible ones. 

\begin{lemma}
	\label{lemma:SL2Ad}
	As $\operatorname{SL}(2,\mathbf{C})$-modules, we have
	\[
	\Lie{sl}(n,\mathbf C) \cong V_{2n - 1} \oplus V_{2n - 3} \cdots \oplus V_3.
	\]
\end{lemma}
\begin{proof}
	Consider the action of $\operatorname{SL}(2,\mathbf C)$  on $\Lie{gl}(n, \mathbf{C})$ obtained by composing the n-dimensional representation $V_n$
	with the adjoint. We have the following isomorphisms of
	$\operatorname{SL}(2,\mathbf C)$-modules:
	\[
	V_n\otimes V_n^*\cong \Lie{gl}(n, \mathbf{C}) \cong \Lie{sl}(n, \mathbf{C}) \oplus \mathbf{C},
	\]
	where the factor $\mathbf{C}$ corresponds to diagonal matrices. 
	The result now follows from the Clebsh-Gordan formula applied to 
	$V_n \otimes V_n^*\cong V_n \otimes V_n$.
\end{proof}

\section{Raghunathan's cohomology vanishing theorem}
\label{sec:vanishing}
The aim of this section is to prove Theorem \ref{thm:l2vanish} stated below.
This theorem is a particular case of a theorem due to Raghunathan \cite{RagAJM}.
Before stating it, let us recall some facts.

The homogeneous manifold $\operatorname{SL}(2, \mathbf C) / \operatorname{SU}(2)$ 
is endowed with a Riemannian structure using the the Killing form on 
$\operatorname{SL}(2, \mathbf C)$ (see section \ref{sec:har_form} for details),
which makes this space isometric to hyperbolic $3$-dimensional space $\mathbf{H}^3$.

Let $\Gamma$ be a discrete torsion-free subgroup of $\operatorname{SL}(2, \mathbf C)$,
and $M = \Gamma \backslash \mathbf{H}^3$ the corresponding complete hyperbolic manifold.
Let $V$ be a finite dimensional representation of $\operatorname{SL}(2, \mathbf C)$,
and $\rho\colon \Gamma \to \operatorname{SL}(V)$ the induced representation.
We can consider the associated \emph{flat} vector bundle over $M$, 
\[ 
E_\rho = \widetilde{M} \times_{\Gamma} V.
\]

The space of $E_\rho$-valued differential forms on $M$ will be denoted by $\Omega^*(M;E_\rho)$. 
A $\operatorname{SU}(2)$-invariant hermitian product on $V$ yields a well defined hermitian metric 
on $E_\rho$, and hence on $\Omega^*(M; E_\rho)$. 
In particular, it makes sense to talk about $L^2$-forms of $\Omega^*(M; E_\rho)$
as those which are square summable.

\begin{theorem}[\cite{RagAJM}]
\label{thm:l2vanish}
Let $\Gamma$ be a discrete torsion-free subgroup of $\operatorname{SL}(2, \mathbf C)$.
Let $V$ be an irreducible finite dimensional complex representation of $\operatorname{SL}(2, \mathbf C)$, 
and $\rho \colon \Gamma \to \operatorname{SL}(V)$ the induced representation.
Then, for $p=1,2$, every closed $L^2$-form in $\Omega^p( \Gamma \backslash \mathbf{H}^3; E_\rho)$ 
 is exact.
\end{theorem}

As an immediate corollary of Theorem \ref{thm:l2vanish} we get a particular 
case of Raghunathan's cohomology vanishing theorem.

\begin{corollary}[\cite{RagAJM}]
\label{cor:closed_l2vanish}
	Let $M$ be a \emph{closed} hyperbolic three-manifold. If $V$ is
	an irreducible finite dimensional complex representation
	of $\operatorname{SL}(2, \mathbf C)$, then 
	\[
	H^1(M; E_\rho) = 0.
	\]
\end{corollary}

\begin{rem}
Raghunathan's theorem  applies to lattices of a semisimple Lie group $G$,
and a broader family of representations, see \cite{RagAJM}. 
\end{rem}

From Theorem \ref{thm:l2vanish} we can easily deduce Theorem \ref{thm:inclusion}.
\begin{proof}[Proof of Theorem~\ref{thm:inclusion}]
	We have $M = \Gamma \backslash \mathbf{H}^3$ for some discrete torsion-free 
	subgroup $\Gamma$ of $\operatorname{SL}(2, \mathbf C)$.
	If $M$ is compact then the result is clear from Theorem \ref{thm:l2vanish}, so
	we can assume $M$ is noncompact.
	The space $H^p(\overline{M},\partial\overline{M}; E_\rho)$ can be identified
	with the cohomology group of compactly supported $E_\rho$-valued $p$-forms 
	on $M$; hence, an element $[\alpha] \in H^p(\overline{M},\partial\overline{M}; E_\rho)$
	is represented by a closed compactly supported form $\alpha$ on $M$.
	Therefore, Theorem \ref{thm:l2vanish} implies that for $p=1,2$
	the image of $[\alpha]$ under the map 
	$H^p(\overline{M},\partial \overline{M}; E_\rho )\rightarrow 
	H^p(M; E_\rho )$ induced by the inclusion is zero. 
	
	The theorem now follows from the long exact sequence of the pair,
	and Poincar\'e duality. Indeed 
	the long exact sequence of the pair $(\overline M, \partial\overline M)$
	gives short exact sequences
	\begin{eqnarray*}
		0\to H^1(M;E_{\rho_n})\to H^1(\partial\overline M;E_{\rho_n})\to H^2( \overline M, \partial{\overline M} ; E_{\rho_n})\to 0,\\
		0\to H^2(M;E_{\rho_n})\to H^2(\partial\overline M;E_{\rho_n})\to H^3( \overline M, \partial{\overline M} ; E_{\rho_n}).
	\end{eqnarray*}
	By Poincar\'e duality we have 
	$\dim  H^1(M; E_{\rho_n}) = \dim H^2( {\overline M}, \partial \overline{M} ; E_{\rho_n})$,
	and $\dim H^3( {\overline M}, \partial \overline{M} ; E_{\rho_n}) = \dim  H^0(M; E_{\rho_n}) = 0$, by Lemma~\ref{lemma:inv}.
\end{proof}

Raghunathan's original proof of the theorem, a particular case of which is Theorem \ref{thm:l2vanish},
uses two results as starting points. The first one the following 
theorem due to Andreotti and Vesentini \cite{AndreottiVesentini}.
Although the original theorem is for complex manifolds, there is 
an adaptation of Garland \cite[Thm.~3.22]{Garland} to the real case.
\begin{theorem}[Andreotti-Vesentini \cite{AndreottiVesentini}, Garland \cite{Garland}]
	\label{thm:AVG}
	Suppose that $M$ is complete. Assume that there exists $c>0$ such that for every
	$\alpha \in \Omega^p(M; E)$ with compact support 
	$(\Delta \alpha, \alpha) \geq c(\alpha, \alpha),$
	where $(\,,)$ denotes the hermitian (or inner) product on the space of $E$-valued forms.
	Then every square-integrable closed $p$-form is exact.
\end{theorem}

The second point is the work of Matsushima-Murakami concerning the theory of
harmonic forms in a locally symmetric manifold \cite{MatMur}.
One of the goals of that work consists in proving a 
Weitzenb\"ock formula for the Laplacian. 
Using that formula, the strong-positivity hypothesis of the 
Laplacian required in Theorem~\ref{thm:AVG} can be proved by 
establishing the positivity of a certain linear operator 
defined on a finite dimensional space, see Subsection \ref{sec:har_form}.
Although this is an important conceptual reduction,
it remains to prove the positivity of that operator.
Raghunathan was able to prove it for a large family of 
locally symmetric manifolds and representations, see \cite{RagAJM}.

The rest of this section is divided into two
parts. The first one is a review the work of Matsushima and Murakami 
concerning the Laplacian of a locally symmetric manifold.
The material presented here is almost entirely based on 
Matsushima-Murakami \cite{MatMur}, and Raghunathan's book \cite{Rag}.
Although it does not bring in a new conceptual approach,
seeking completeness,
we hope the exposition given here will be more accessible
to the non-expert. Using this material, 
we give a simple proof of Theorem \ref{thm:l2vanish}
 in Subsection~\ref{section:proof_rag}.

\subsection{Review of harmonic forms on a locally symmetric manifold }
\label{sec:har_form}

Let $G$ be a connected semi-simple Lie group and $K < G$ a maximal compact subgroup of $G$.
The respective Lie algebras are denoted by $\Lie{g}$ and $\Lie{k}$, with the convention that
they are the Lie algebras of left invariant vector fields on $G$ and $K$, respectively.

Let $B$ denote the Killing form of $\Lie{g}$. We recall that it is defined by 
$$B(V, W) = \operatorname{tr}(\operatorname{ad}_V\circ\operatorname{ad}_W),
$$
 for $V,W\in \Lie{g}$.
Cartan's criterion implies that $B$ is nondegenerate if, and only if, $\Lie{g}$ is semisimple.
In that case, we have a canonical decomposition $ \Lie{g} = \Lie{m} \oplus \Lie{k} $, 
where $\Lie{m}$ is the orthogonal complement to $\Lie{k}$ respect
to $B$. This decomposition satisfies the following properties: $B$ is negatively defined on $\Lie{k}$; 
$B$ positively defined on $\Lie{m}$; $[ \Lie{k}, \Lie{m} ] \subset \Lie{m} $; and
$[ \Lie{m}, \Lie{m} ] \subset \Lie{k} $. 

The Killing form defines a pseudo-Riemannian metric on $G$, 
which is invariant by the action of $G$ by right translations, and is positively (resp.\; negatively) 
defined on $\Lie{m}$ (resp.\;$\Lie{k}$).
Therefore, the Killing form defines a Riemannian metric on the homogeneous space $X = G/K$. Note that
$G$ acts on the left on $X$ by orientation preserving isometries.  

Let $\Gamma$ be a discrete subgroup of $G$ that acts freely on $X$. Since $\Gamma$ acts by isometries, 
the quotient $M = \Gamma \backslash X$ is a Riemannian manifold. It is said that $M$ is a 
\textit{locally symmetric manifold}.

For our purposes, it will be convenient to regard the universal covering
$X\to M$ as a principal bundle over 
$M$ with structure group $\Gamma$.
We follow the convention that the action of the structure group of a principal bundle is on the right. Hence
we only need to convert the action of $\Gamma$ into a right action (if $g\in \Gamma$, then 
$x\cdot g = g^{-1}\cdot x$, for $x \in X$). We will also regard $X$ as a flat bundle.

Consider the $G$-principal bundle $P = X \times_\Gamma G$ over $M$ (see Appendix~\ref{principalbundles} for notation) endowed with 
the flat connection induced from the trivial connection of the product $X \times G$. We can embed $X$ on $P$ using the 
section $X\rightarrow X\times G$ whose second coordinate is constant and equal to the identity element. We can think of
$X$ as a reduction of the structure group. Obviously, the horizontal leaves of $X$ are also horizontal leaves
of $P$, so the connection on $P$ is reducible to $X$. 

On the other hand, the principal bundle $P$ has a canonical reduction of its structure 
group from $G$ to $K$. In order to get such a reduction, consider the embedding
$i\colon G \hookrightarrow X\times G$ given by $i(g) = (gK, g)$. The image of $G$
by this embedding is invariant by the bundle action of $K$, so it defines 
an embedding $\Gamma\backslash G \hookrightarrow X\times_\Gamma G$, which will be
also denoted by $i$. Therefore, $Q = i( \Gamma\backslash G)$ is a 
reduction of the structure group.

The connection defined on $P$ is not reducible to $Q$, because its horizontal distribution
is not tangent to $Q$ (a curve on $X\times G$ whose second component is constant, gives 
an horizontal curve on $P$; hence, if the horizontal distribution were tangent to $Q$, 
this curve would be contained in $Q$, and this does not happen). 
Nevertheless, since the action of $K$ on $\Lie{g}$ respects the decomposition 
$\Lie{g} = \Lie{m}\oplus\Lie{k}$, we can state the following.

\begin{proposition}
	Let $\eta \in \Omega^1(P;\Lie{g})$ be the connection form of the connection 
	defined on $P$ above. Put $\eta = \eta_\Lie{m} + \eta_\Lie{k}$, where $\eta_\Lie{m}$ 
	and $\eta_\Lie{k}$ are the $\Lie{m}$ and $\Lie{k}$ components of $\eta$ respectively. 
	Then, the restriction of $\eta_\Lie{k}$ to $Q$ is a connection form on $Q$.
\end{proposition}

\begin{obs}
	We can identify $\Lie{g}$ with the space of vector fields
	on $\Gamma \backslash G$ that are projection of left invariant vector fields on $G$. 
	In what follows, we will tacitly do this identification.
\end{obs}

Let $\omega \in \Omega^1(\Gamma \backslash G; \Lie{g})$ be the left Maurer-Cartan form of $G$.
It is easily checked that $i^*(\eta) = \omega$. Hence, if we decompose
$\omega = \omega_\Lie{m} + \omega_\Lie{k}$ into the $\Lie{m}$-component and the $\Lie{k}$-component,
$\omega_\Lie{k}$ is the connection form of the connection defined on $\Gamma \backslash G$, and
the horizontal distribution is given by $\Lie{m}$.

Consider a finite linear representation $\rho\colon G\longrightarrow \text{Aut}(V),$ and 
the associated vector bundle $E = X \times_\Gamma V$ (note that $E$ is canonically 
identified with $P \times_G V$ and $Q \times_K V$). 

The flat connection on $P$ defines an exterior covariant differential $d_\rho$ 
on the space $\Omega^*(M; E)$. Via the canonical isomorphism between 
$\Omega^*_{\operatorname{Hor}}(\Gamma\backslash G;V)^K$ and 
$\Omega^*(M; E)$, we can transfer the operator $d_\rho$ to an operator $D_\rho$, 
in such a way that this isomorphism is a chain complex isomorphism. If we denote by $D$
the exterior covariant differential defined by the connection $\omega_\Lie{k}$ on $Q$, 
then the relation between $D$ and $D_\rho$ is given by the following proposition.

\begin{proposition}
	\label{canonical_iso}
	Let $\alpha$ be a form in $\Omega^r_{\operatorname{Hor}}(\Gamma\backslash G;V)^K$.
	We have the following decomposition 
	$$D_\rho\alpha = D \alpha + T \alpha,$$ where $T\alpha = \rho(\omega_{m})\wedge\alpha$.
\end{proposition}
\begin{proof}
	On $P$ the differential covariant is given by $d\alpha + \rho(\eta)\wedge\alpha$ 
	(see Proposition \ref{prop:dif_cov}). Hence, if we transfer it to 
	$Q$ via $i$, we get $D_\rho \alpha = d\alpha + \rho(i^*\eta)\wedge\alpha,$ and the proposition
	follows from the fact that $i^*\eta = \omega$.
\end{proof}

Let's fix an orientation on $\Lie{k}$ and $\Lie{m}$, and take an orthonormal basis for 
$\Lie{g}$, $(X_1,\dots,X_n,Y_1,\dots,Y_m)$, such that 
$(X_1,\dots,X_n)$ and $(Y_1,\dots,Y_m)$ are positively oriented orthonormal bases for $\Lie{k}$ and $\Lie{m}$, respectively. 
Here, orthonormality means that  
$$ B(X_i,X_j) = -\delta_{ij} \quad B(Y_i,Y_j) = \delta_{ij}, \quad B(X_i,Y_j) = 0. $$

\begin{notation}
	We will follow the following conventions. Let $V$ be a finite dimensional vector space.
	If $e_1,\dots,e_n$ is a basis for $V$, then its dual basis will be denoted by $e^1,\dots,e^n \in  V^*$, with
	$e^i(e_j) = \delta_{ij}$. 
	If $A \in \bigotimes^r V^*$ is an $r$-times covariant tensor, then 
	its components relative to the basis defined by $e^1,\dots, e^n$ will be denoted
	by $A_{i_1,\dots,i_r}$. 
	Concerning the exterior product on $\bigwedge^* V^*$, we will follow the
	convention such that $e^1\wedge\cdots\wedge e^n$ is the determinant. We will also use Einstein notation. 
	Hence, given $\alpha \in \bigwedge^r V^*$, we have
	$\alpha = \alpha_{i_1,\dots,i_r} e^{i_1}\otimes\cdots\otimes e^{i_r},$ where $\alpha_{i_1,\dots,i_r}$ are scalars
	satisfying $\alpha_{i_{\sigma(1)},\dots,i_{\sigma(r)}} = \operatorname{sgn}(\sigma)\alpha_{i_1,\dots,i_r},$ for 
	any permutation $\sigma \in \Sigma_r$. Then we also have
	$$\alpha = \sum_{1\leq i_1<\cdots<i_r\leq r} \alpha_{i_1,\dots,i_r} e^{i_1}\wedge\cdots\wedge e^{i_r}
	         = \frac{1}{r!}\alpha_{i_1,\dots,i_r} e^{i_1}\wedge\cdots\wedge e^{i_r}.$$
\end{notation}

From now on, all the tensors will be written in the basis of $\Lie{g}$ 
given by $\{X_1,\dots,X_n,Y_1,\dots,Y_m\}$.

\begin{proposition}
	\label{prop:frame_expr}
	For $\alpha \in \Omega^r_{\operatorname{Hor}}(\Gamma\backslash G;V)^K$, 
	the operators $D$ and $T$ are given by the following equations.
	\begin{eqnarray}
		(D\alpha)_{i_1,\dots,i_{r +1}} & = & \sum_{k=1}^{r+1} (-1)^{k+1} Y_{i_k} \alpha_{i_1,\dots,\widehat{i_k},\dots,i_{r+1}} 
		\label{eqn_D}\\
		(T\alpha)_{i_1,\dots,i_{r +1}} & = & \sum_{k=1}^{r+1} (-1)^{k+1} \rho(Y_{i_k}) \alpha_{i_1,\dots,\widehat{i_k},\dots,i_{r+1}}
		\label{eqn_T}.
	\end{eqnarray}
\end{proposition}
\begin{proof}
	Put $\alpha = \frac{1}{r!}\alpha_{i_1,\dots,i_r} Y^{i_1}\wedge\cdots\wedge Y^{i_r}.$
	By definition, $D\alpha$ is the horizontal component of $d\alpha$. It is immediate that $dY^k$ has no 
	horizontal component: $dY^k(Y_i, Y_j) = Y^k([Y_i,Y_j]) = 0$. Hence, 
	$D\alpha = \frac{1}{r!} Y_j \alpha_{i_1,\dots,i_r} \otimes Y^j\wedge Y^{i_1}\wedge\cdots\wedge Y^{i_r}$. 
	Rearranging the indices we get equation \ref{eqn_D}. The other equation follows immediately from the
	definition of $T$.
\end{proof}

Let us define the forms $\Omega_K = X^1\wedge\cdots\wedge X^n$ and $\Omega_M = Y^1\wedge\cdots\wedge Y^m$.
It is clear that these forms are independent of the orthonormal bases chosen. 
Hence, $\Omega_K$ and $\Omega_M$ are well defined forms on $\Gamma\backslash G$. Note that
$\Omega_K$ is vertical and $\Omega_M$ is horizontal, and both are right $K$-invariant 
(it is a consequence of the fact the right action of $K$ on $\Lie{g}$ leaves both the Killing form
and the decomposition $\Lie{g} = \Lie{k}\oplus\Lie{m}$ invariant). Observe that $\Omega_M$ defines 
a volume form on $M$, which is compatible with the metric structure of $M$.

Next we want to define an inner product on the fibers of $E$. In order to do that, fix a $K$-invariant inner
product $\langle\;,\rangle_V$ on $V$, and use it to define a metric on the fibers of $E = Q\times_K V$. 
Then define an inner product on 
$\Omega^*(M; E)$ as usual: if $\alpha,\beta \in \Omega^*(M; E)$ then
$$(\alpha,\beta) = \int_M \langle\alpha(x),\beta(x)\rangle_x \Omega_M,$$
where $\langle\,,\rangle_x$ is the inner product defined on the fiber $E_x$, and $\Omega_M$ is interpreted as
a form on $M$. On the other hand, we can define the inner product of two forms 
$\tilde{\alpha},\tilde{\beta}\in \Omega_{\operatorname{Hor} }^r(\Gamma \backslash G;V)^K$ by
$$
	(\tilde{\alpha},\tilde{\beta}) =
	 \frac{1}{\mu(K)} \int_{\Gamma \backslash G}\langle \tilde{\alpha}(u), \tilde{\beta}(u) \rangle_u \Omega_K\wedge\Omega_M,
$$
where $\langle\,,\rangle_u$ is the inner product on $\bigwedge^r H^*\otimes V$ induced by the the Killing form, and the 
inner product on $V$, and $\mu(K) = \int_K \Omega_K$ the volume of $K$. 
Proposition \ref{G_bundle:int}, gives the relation between these two products.

\begin{proposition}
	The canonical isomorphism between $\Omega_{\operatorname{Hor} }^*(\Gamma \backslash G;V)^K$
	and $\Omega^*(M; E)$ is an isometry.
\end{proposition}

Using the Hodge dual operator on the horizontal bundle 
$$* \colon \Omega_{\operatorname{Hor} }^r(\Gamma \backslash G;V)^K\rightarrow 
\Omega_{\operatorname{Hor} }^{m-r}(\Gamma \backslash G;V)^K,$$ 
we can give a characterization of the formal adjoint
of the operators $D$ and $T$.
\begin{proposition}
	Let $\alpha\in\Omega_{\operatorname{Hor} }^r(\Gamma \backslash G;V)^K$ with compact support. Then,
	\begin{eqnarray}
		D^*\alpha & = & (-1)^{r}*^{-1} D *\alpha \label{eqn:D*} \\
		T^*\alpha & = & (-1)^{r - 1 }*^{-1}\overline{\rho(\omega)}^t\wedge(*\alpha) \label{eqn:T*} 
	\end{eqnarray}
\end{proposition}

\begin{proof}
	We want to use Proposition \ref{G_bundle:adjoint}. We claim that 
	$$\int_P D\alpha\wedge\beta\wedge\Omega_K = (-1)^r\int_P \alpha\wedge D \beta\wedge\Omega_K,$$
	for $\alpha$ and $\beta$ forms of $\Omega_{\operatorname{Hor} }^*(\Gamma \backslash G;V)$ with compact support of degree $r-1$ 
	and $m-r $ respectively. Indeed, since $D\alpha$ is the horizontal component of $d\alpha$, we have 
	$D\alpha\wedge\Omega_K = d\alpha\wedge\Omega_K$. Then,
	$$ d(\alpha\wedge\beta\wedge\Omega_K) = d\alpha\wedge\beta\wedge\Omega_K + (-1)^{r-1}\alpha\wedge d\beta\wedge\Omega_K,$$
	for $\Omega_K$ being closed. Therefore, by Stokes' Theorem we get the equality we wanted to prove. Now, Proposition
	\ref{G_bundle:adjoint} gives Formula~(\ref{eqn:D*}).

	Now, let us prove (\ref{eqn:T*}). By Proposition \ref{G_bundle:adjoint}, it suffices
	to prove that $$(\rho(\omega)\wedge\alpha)\wedge\beta = (-1)^{r-1} \alpha\wedge(\rho(\omega)^*\wedge\beta).$$
	If we take an orthonormal basis for $V$, then $\alpha$ and $\beta$ are column vectors of forms of degree $r-1$ and $m - r$ respectively, 
	and $\rho(\omega)$ a matrix of one forms. Hence, in this basis $(\rho(\omega)\wedge\alpha)\wedge\beta$ is
	$(\rho(\omega)\alpha)^t\bar{\beta} $, but $(\rho(\omega)\alpha)^t\beta = (-1)^{r-1}\alpha^t \rho(\omega)^t \overline{\beta}$, as we wanted to 
	prove.

\end{proof}
A similar proof of Proposition \ref{prop:frame_expr}, using the formulae found in the previous proposition, gives the following.
\begin{proposition}
	\label{prop:frame_expr_st}
	For $\alpha \in \Omega^r_{\operatorname{Hor}}(\Gamma\backslash G;V)^K$ with compact support, 
	the operators $D^*$ and $T^*$ are given by the following equations.
	\begin{eqnarray}
		(D^*\alpha)_{i_1,\dots,i_{r-1}} & = & \sum_{k=1}^{m} - Y_{k} \alpha_{k,i_1,\dots,i_{r-1}} \label{eqn_D_st} \\
		(T^*\alpha)_{i_1,\dots,i_{r-1}} & = & \sum_{k=1}^{m} \rho(Y_k) \alpha_{k,i_1,\dots,i_{r-1}} \label{eqn_T_st}.
	\end{eqnarray}
\end{proposition}

\begin{lemma}
	\label{lemma:S=0}
	If the inner product on $V$ is symmetric respect to the action of $\Lie{m}$, then
	the operator $S = TD^* + T^*D + DT^* + D^* T$ is zero for every form with compact support. 	
\end{lemma}

Before proving the lemma, we need the following result.
\begin{lemma}
	\label{lemma:van_int}
	For every function $f$ with compact support, and $Y\in \Lie{g}$, 
	$$
	\int_{\Gamma\backslash G} (Yf) \Omega_M\wedge\Omega_K = 0.
	$$ 
\end{lemma}
\begin{proof}
	Since $Y$ is an infinitesimal isometry we have $L_Y( f \Omega_M\wedge\Omega_K) = (Yf)\Omega_M\wedge\Omega_K$.
	On the other hand, the formula $L_Y = i_Y\circ d + d\circ i_Y$ gives $L_Y( f \Omega_M\wedge\Omega_K) = d( i_Y f \Omega_M\wedge\Omega_K),$
	and Stokes' Theorem implies 
	   $$ 0 = \int_{\Gamma\backslash G} L_Y(f \Omega_M\wedge\Omega_K) = \int_{\Gamma\backslash G}(Yf)\Omega_M\wedge\Omega_K,$$ 
	as we wanted to prove.
\end{proof}

\begin{proof}[Proof of Lemma~\ref{lemma:S=0}]
	Since $S$ is a self-adjoint operator, $S = 0$ if, and only if, 
	$(S\alpha,\alpha) = 0$ for every $\alpha$ with compact support. 
	Let's take $\alpha \in \Omega_{\operatorname{Hor} }^*(\Gamma \backslash G;V)^K$ 
	with compact support. We must show that	
	$$(S\alpha,\alpha) =(D\alpha, T\alpha) + (T\alpha, D\alpha) + (D^*\alpha, T^*\alpha) + (T^*\alpha, D^*\alpha)  = 0.$$
	Observe that it suffices to prove that $(D\alpha, T\alpha) + (D^*\alpha, T^*\alpha)  = 0$.
	Moreover, using the $\Lie{m}$-symmetry of the inner product and the fact that the Hodge $*$ operator 
	is an isometry,  we must prove $(D\alpha, T\alpha) + (D(*\alpha), T(*\alpha))  = 0.$ 
	Let's compute $(D\alpha, T\alpha)$. Put $\alpha = \alpha_{i_1,\dots,i_r}\otimes Y^{i_1}\wedge\cdots\wedge Y^{i_r}.$ 
	If we use the expression of $D$ and $T$ given in Proposition \ref{prop:frame_expr}, we see that $(D\alpha, T\alpha)$ 
	is the sum of terms of the form 
	$$(-1)^{i+j}\int_{\Gamma\backslash G} \langle Y_{i_j} \alpha_{i_1,\dots,\widehat{i_j},\dots,i_{r+1}},
	\rho(Y_{i_k})\alpha_{i_1,\dots,\widehat{i_k},\dots, i_{r+1}} \rangle_V d\mu_G.$$
	It is convenient to group the summands according to whether the avoided sub-indices
	$\widehat{i_j}$ and $\widehat{i_k}$ are equal or not.
	Therefore, one term is a sum of factors of the form
	$$ \int_{\Gamma\backslash G} \langle Y_j \alpha_{i_1,\dots,i_r}, 
	\rho(Y_j)\alpha_{i_1,\dots, i_r} \rangle_V d\mu_G, \quad j\notin \{i_1,\dots, i_r\},$$
	and the rest is a sum of terms of the form
	\begin{equation}
		\label{eq_DT}
		(-1)^{j + k} \int_{\Gamma\backslash G} \langle Y_{i_j} \alpha_{i_1,\dots,\widehat{i_j},\dots, i_k,\dots, i_r},
		\rho(Y_{i_k})\alpha_{i_1,\dots,i_j,\dots, \widehat{i_k},\dots, i_r}\rangle d\mu_G,
	\end{equation}
	with $i_j \neq i_k$. We can apply this formula to $*\alpha$ to compute $(D(*\alpha), T(*\alpha))$. The formula we get
	is just the above formula with the range of the indices changed by their complementary; that is, one one hand we get terms
	of the form 
	$$ \int_{\Gamma\backslash G} \langle Y_j \alpha_{i_1,\dots,i_r}, 
	\rho(Y_j)\alpha_{i_1,\dots, i_r} \rangle_V d\mu_G,\quad j\in \{i_1,\dots, i_r\},$$
	and on the other hand terms of the form
	$$ (-1)^{j + k} \int_{\Gamma\backslash G} \langle Y_{i_k} \alpha_{i_1,\dots,\widehat{i_j},\dots, i_k,\dots, i_r},
	\rho(Y_{i_j})\alpha_{i_1,\dots,i_j,\dots, \widehat{i_k},\dots, i_r}\rangle d\mu_G,$$
	for $i_j\neq i_k$. By Lemma \ref{lemma:van_int}, this last term is the opposite of \ref{eq_DT}.
	Hence, it suffices to prove that for every $Y\in \Lie{m}$, and $f \in \mathcal{C}(\Gamma\backslash G; V)$, we have 
	$$ 
	\int_{\Gamma\backslash G} \langle Y f, \rho(Y) f \rangle_V d\mu_G=0.
	$$
	But it is also an immediate consequence of Lemma \ref{lemma:van_int} and the symmetry of $\rho(Y)$.
	The lemma now follows from the fact that $(D^*\alpha, T^*\alpha) = (D(*\alpha), T(*\alpha))$.
	
\end{proof}

\begin{corollary}[Matsushima-Murakami formula]
	\label{coro:MMformula}
	Assume the inner product on $V$ is symmetric respect to the action of $\Lie{m}$. Then
	$$
	\Delta_\rho = \Delta + H_\rho,
	$$
	where $ \Delta = D D^* + D^*D$, and
	$H_\rho = TT^* + T^*T$. 
\end{corollary}

\begin{proof}
	We have $\Delta_\rho = D_\rho D_\rho^*+D_\rho^* D_\rho = \Delta + H_\rho + S$, and
	Lemma~\ref{lemma:S=0}.
\end{proof}

Let's denote by $\mathbf{T}$, $\mathbf{T}^*$, $\mathbf{H}_\rho$ the restriction to $V\otimes\bigwedge^p\Lie{m}^*$
of $T$, $T^*$ and $H_\rho$ respectively. Since $T$ is an operator of degree zero, 
essentially all information of $T$, $T^*$ and $H_\rho$ is contained in $\mathbf{T}$, $\mathbf{T}^*$, $\mathbf{H}_\rho$.
In particular, $H_\rho$ is positive definite if and only  $\mathbf{H}_\rho$ is so.

\begin{proposition}
\label{prop:Def_H}
	Let $\alpha \in V\otimes\bigwedge^p\Lie{m}^*$. Then we have, 
	\begin{eqnarray*}
		(\mathbf H_p \alpha )_{i_1,\dots,i_r} & = & \sum_{j= 1}^m \rho(Y_j)^2 \alpha_{i_1,\dots,i_r}  + 
		\sum_{k=1}^r \sum_{j= 1}^m (-1)^{k+1}\rho([Y_{i_k},Y_j])\alpha_{j,i_1,\dots,\widehat{i_k},\dots,i_r} 
	\end{eqnarray*}
\end{proposition}

\begin{proof}
	Put $\beta_{i_1,\dots,i_{r+1} } = (T\alpha)_{i_1,\dots,i_{r+1} } $ and  
	$\gamma_{i_1,\dots,i_{r-1}} = (T^*\alpha)_{i_1,\dots, i_{r-1}} $. Then, on one hand we have
	\begin{eqnarray*}
		(TT^*\alpha)_{i_1,\dots,i_r} & = & 
		\sum_{k = 1}^r (-1)^{k+1}\rho(Y_{i_k}) \gamma_{i_1,\dots,\widehat{i_k},\dots,i_r} \\
		& = & \sum_{k = 1}^r (-1)^{k+1} \rho(Y_{i_k}) \sum_{j= 1}^m \rho(Y_j) \alpha_{j,i_1,\dots,\widehat{i_k},\dots,i_r}. 
	\end{eqnarray*}
	and on the other hand,

	\begin{eqnarray*}
		(T^*T \alpha)_{i_1,\dots,i_r}&=& \sum_{j = 1}^m \rho(Y_j) \beta_{j,i_1,\dots,i_r}\\
		& = & \sum_{j = 1}^m \rho(Y_j)\big( \rho(Y_j)\alpha_{i_1,\dots,i_r} + \sum_{k=1}^r (-1)^k \rho(Y_{i_k})\alpha_{j,i_1,\dots,\widehat{i_k},\dots,i_r} \big). 
	\end{eqnarray*}
	And the proposition follows.
\end{proof}

\subsection{Proof of Theorem~\ref{thm:l2vanish}}
\label{section:proof_rag}

We want to apply the criterion of Andreotti-Vesentini of  Theorem~\ref{thm:AVG}.
For this purpose, we will use  Matsushima-Murakami's formula (Corollary~\ref{coro:MMformula})
for the representation of $\operatorname{SL}(2,\mathbf C)$. 
Since for every compactly supported 1-form $\alpha$
\[
(\Delta(\alpha),\alpha)=(D(\alpha),D(\alpha))+(D^*(\alpha),D^*(\alpha))\geq 0,
\]
using Corollary~\ref{coro:MMformula}, the criterion of
Theorem~\ref{thm:AVG} reduces to show 
that $(H_{\rho}(\alpha),\alpha)\geq c(\alpha,\alpha)$ for some 
uniform  $c>0$ and every compactly supported 1-form $\alpha$.

Notice that since 
the linear operator $H_{\rho}$ on 1-forms is induced from
a linear operator $\mathbf{H}_{\rho}$ on 
$V\otimes \mathfrak m^*$, if $\mathbf{H}_{\rho}$ is positive definite,
then there is a positive constant $c$ so that 
$(H_{\rho}(\alpha),\alpha)\geq c(\alpha,\alpha)$ holds for every 
compactly supported one form $\alpha$.
The proof will follow from Lemma~\ref{lemma:Hpositivedef}.

In order to apply Matsushima-Murakami's formula 
to the representations of $\operatorname{SL}(2,\mathbf C)$, 
first we need to choose an 
orthonormal basis for $\Lie{su}(2)$ respect to the Killing form 
(in fact, respect to a constant multiple of it). Let's define
\[
X_1 = \left ( \begin{array}{cc}
	\mathbf{i} &  0 \\
	0 & -\mathbf{i}
	\end{array} 
	\right),
X_2 = \left ( \begin{array}{cc}
	0  & 1\\
	-1 & 0
\end{array} \right),
X_3 = \left ( \begin{array}{cc}
	0 & \mathbf{i}\\
	\mathbf{i} & 0
\end{array} \right).
\]
Then $(X_1, X_2, X_3)$ is an orthonormal basis for $\Lie{su}(2)$.
The orthogonal complement to 
$\Lie{su}(2)$ with respect to 
the Killing form is given 
by $Y_k = \mathbf{i} X_k$, for $k = 1, 2, 3$. On the other hand, we have 
$[X_i, X_{i+1}] = 2 X_{i+2}$, for $i = 1,2,3$, 
where the indices are taken modulo $3$.

\begin{lemma}
	\label{lemma:Hpositivedef}
	Let   $\rho \colon \Lie{sl}(2, \mathbf{C})\rightarrow \operatorname{End}(V)$     a \emph{complex} finite dimensional irreducible representation,
	$\dim(V)\geq 2$.
	Then the operator $\mathbf{H}_{\rho}$ is positively defined on degree $1$ and $2$.
\end{lemma}

\begin{proof}
	Since $\mathbf{H}_{\rho}=\mathbf{T}_{\rho} \mathbf{T}^*_{\rho}+\mathbf{T}^*_{\rho} \mathbf{T}_{\rho}$, to show that 
	$\mathbf{H}_{\rho}$ is positive definite is equivalent to show that its kernel is trivial.
	Let $\alpha\in V\otimes \Lie{m}^*$. We have  
	$\alpha = \sum_{i = 1}^3\alpha_{i}\otimes Y^i$, with $\alpha_i \in V$. 
	Assume $\mathbf{H}_\rho \alpha = 0$. Then $\mathbf{T}_\rho\alpha = 0$ must vanish too, 
	and from Proposition~\ref{prop:frame_expr}~(\ref{eqn_T}) we obtain 
	\begin{equation}
		0 = (\mathbf{T}_\rho\alpha)(Y_i,Y_j)=\rho(Y_i)\alpha_j-\rho(Y_j)\alpha_i,\quad i,j = 1,2,3.
		\label{rel_T}
	\end{equation}
	Proposition~\ref{prop:Def_H} yields
	$$(\mathbf{H}_\rho \alpha)(Y_j)=\sum_{k=1}^3 \big(\rho(Y_k)^2 \alpha_j+\rho([Y_j,Y_k])\alpha_k\big).$$
	Taking the indices modulo $3$, and using the Lie algebra relations, we get
	\begin{eqnarray*}
		\sum_{k=1}^3 \rho([Y_j,Y_k])\alpha_k& = 
		&\rho([Y_j,Y_{j+1}])\alpha_{j+1}+\rho([Y_j,Y_{j+2}])\alpha_{j+2}\\
		& = & 2\big(\rho( - X_{j+2} )\alpha_{j+1} + \rho( X_{j+1} ) \alpha_{j+2} \big ) \\
		& = & 2\mathbf{i}\big(\rho( Y_{j+2} )\alpha_{j+1} - \rho( Y_{j+1} ) \alpha_{j+2} \big ).
	\end{eqnarray*}
	Notice that in the last equality we have used the complex structure. 
	Hence, using \ref{rel_T}, we get $(\mathbf{H}_\rho \alpha)(Y_j)=\sum_{k=1}^3 \rho(Y_k)^2 \alpha_j$, and then
	\begin{eqnarray*}
		0 = \langle \mathbf{H}_\rho\alpha,\alpha\rangle& = & \sum_{j=1}^3
		\big\langle\sum_{k=1}^3 \rho(Y_k)^2 \alpha_j,\alpha_j\big\rangle \\
		& = &\sum_{j,k=1}^3 \big \langle \rho(Y_k)\alpha_j,\rho(Y_k)\alpha_j\big\rangle,
	\end{eqnarray*}
	that implies $\rho(Y_j)\alpha_k = 0$  for $j,k=1,2,3$.
	Hence, for a fixed $k$, we have $\rho(Z) \alpha_k = 0$ for every $Z\in\Lie{sl}(2,\mathbf{C})$. 
	Since we are assuming that $\rho$ is irreducible and nontrivial, we get $\alpha_k = 0$ for all $k$.
	It proves the the lemma in degree $1$. 
	Since $\Lie{m}^*\cong \bigwedge^2\Lie{m}^*$, the same proof holds true in degree $2$.
\end{proof}

\section{Cohomology of the ends and lifts of the holonomy}
\label{sec:cohom_ends_lifts}

Assume that $M$ is a noncompact, nonelementary, orientable hyperbolic manifold with finite topology,
in particular it is the interior of a compact manifold with boundary $\partial\overline M$.
The aim of this section is to analyse the cohomology groups of 
$H^*( \partial \overline M, E_{\rho_n} )$. This will be done in Subsection~\ref{sec:cohom_ends}.
When the ends of the manifold are cusps, this cohomology happens to be related to the lift of the holonomy, 
that we study in Subsection~\ref{sec:lift}.
Finally, this is used to prove Theorem~\ref{thm:H1cusp}.

\subsection{Cohomology of the ends}
\label{sec:cohom_ends}

\begin{definition}
	Let $G$ be a group acting on a vector space $V$. 
	The \emph{subspace of invariants} of $V$, denoted by $V^G$, 
	is the subspace consisting of elements of $V$ that are fixed by $G$.
	That is, 
	\[
	V^G=\{ v \in V \mid g \cdot v = v, \textrm{ for all } g \in G \}.
	\]
\end{definition}

\begin{lemma}
\label{lemma:H_F}
	Let $F$ be a connected component of $\partial \overline M $. 
	For every $n>1$ we have,
	\begin{eqnarray*}
		\dim_\mathbf{C} H^0( F; E_{\rho_n} ) & = &  \dim_\mathbf{C} V_n ^{\pi_1(F)} \\
		\dim_\mathbf{C} H^1( F; E_{\rho_n} ) & = & 2\dim_\mathbf{C} V_n ^{\pi_1(F)} - n \chi(F), \\ 
		\dim_\mathbf{C} H^2( F; E_{\rho_n} ) & = &  \dim_\mathbf{C} V_n ^{\pi_1(F)}.
	\end{eqnarray*}	
\end{lemma}
\begin{proof}
	Since $F$ is a $K(\pi_1(F), 1)$ space, $H^0( F; E_{\rho_n} ) = H^0( \pi_1(F); E_{\rho_n} ),$
	and this is identified with $V_n ^{\pi_1(F)}$. 
	It proves the first equality. The third one follows from Poincar\'e duality,
	and the second one from an Euler characteristic argument. 
\end{proof}

Therefore, all the cohomological information comes from the 
subspace of invariants $V_n ^{\pi_1(F)}$. We distinguish two cases according to whether
$F$ has genus $g\geq 2$, or $F$ is a torus. In order to analyse the 
case when $F$ is a torus, we make the following definition.
If we have a torus $T^2 \subset \partial \overline M$, then
the holonomy  maps $\pi_1 (T^2)$   to a parabolic subgroup;
hence, up to conjugation every element in $\pi_1(T^2)$ is mapped by a lift 
of the holonomy representation to
\[
\pm\begin{pmatrix}
     1 & * \\
    0 & 1
    \end{pmatrix}.
\]

\begin{definition}
Let us fix a lift to $\operatorname{SL}(2,\mathbf C)$ of the holonomy representation.
We say that this lift is \emph{positive} on $\pi_1(T^2)$ 
if every element of $\pi_1(T^2)$ has trace $+2$.
\end{definition}

\begin{proposition}
	\label{prop:invariant}
	Let $F$ a connected component of $\partial \overline M$, and $n > 1$.
	If $F$ has genus  $g\geq 2 $, then $V_n^{\pi_1(F_g)} = 0$.

	If $F$ is a torus $T^2$, then we have the following cases,
	\[
	V_n^{\pi_1(T^2)} = 
	\left\{ 
		\begin{array}{ll}
			0 &  \textrm{ for } n \textrm{ even and a nonpositive lift;} \\
			\mathbf C &  \textrm{ for } n \textrm{ even and a positive lift;} \\
			\mathbf C &  \textrm{ for } n \textrm{ odd.}
		\end{array} 
	\right..
	\]
\end{proposition}

Before proving it, we need the following lemmas. The first one can be found in standard references about Kleinian groups (cf.\ \cite{Kapovich}):

\begin{lemma}
\label{lemma:elementary}
Let $M$ be a hyperbolic three manifold. Then the following are equivalent:
\begin{itemize}
 \item[--] $M$ is elementary (its holonomy is reducible in $\operatorname{PSL}(2,\mathbf C)$).
 \item[--] $\pi_1(M)$ is abelian.
 \item[--] $M$ is   homeomorphic to either the product of the  plane with a circle, $\mathbf{R}^2\times S^1$, or to the product of 
 a 2-torus with a line, $S^1\times S^1\times\mathbf R$.
\end{itemize}
\end{lemma}

\begin{lemma}
 \label{lemma:bdry}
Let $F$ be a connected component of $\partial \overline M $. If $F$ has genus  $g\geq 2$, 
then  $\widetilde{\operatorname{Hol}}(\pi_1 (F))$ is an \emph{irreducible} subgroup of $\operatorname{SL}(2,\mathbf C)$.
\end{lemma}

\begin{proof}
	When $F$ is $\pi_1$-injective (i.e.\ when $\pi_1(F)$ injects into $\pi_1(M)$) then the holonomy
	restricts to a discrete and faithful representation of $\pi_1(F)$,
	and irreducibility follows because $\pi_1(F)$ is nonabelian.
	Otherwise, when $F$ is not $\pi_1$-injective, according to Bonahon \cite{Bonahon} and  McCullough-Miller \cite{McCM} there are two possibilities: 
	either $M$ is a \emph{handlebody} or $F$ is a boundary component of a  \emph{characteristic compression body} $C\subseteq M$. 
	A handlebody  is the result of attaching one handles to a 3-ball; in particular when $M$ is a handlebody then $\pi_1(F)$ surjects 
	onto $\pi_1(M)$, thus $\operatorname{Hol}(\pi_1(F))= \operatorname{Hol}(\pi_1(M))$ and irreducibility comes from the hypothesis 
	that $M$ is nonelementary. 
	Next, assume that $F$ is the positive boundary of a characteristic compression body $C$, 
	namely  $C\subseteq M$  is a codimension $0$ closed submanifold,
	whose boundary splits as a union $\partial C=\partial_- C\cup \partial _+ C$, so that $\partial_+C=F$, the components of 
	$\partial_- C$ are  $\pi_1$-injective in $M$,
	and $C$ is the result of gluing $1$-handles to $\partial_- C\times [0,1]$ along  $\partial_- C\times \{ 1\}$. 
	In particular $\pi_1(F)$ surjects onto $\pi_1(C)$ and $\operatorname{Hol}(\pi_1 (F))=\operatorname{Hol}(\pi_1 (C))$. Thus, if $F=\partial_+C $ 
	and one of the components of $\partial_- C$
	has genus $\geq 2$, then we are done by the $\pi_1$-injective case. Finally if $F =\partial_+C $ 
	and all components of $\partial_- C$ are tori, since incompressible tori in $M$ are boundary parallel, 
	then the inclusion $C\subseteq M$ is a homotopy 
	equivalence. Thus  $\pi_1(F)$ surjects onto $\pi_1(M)$ and irreducibility follows again because $M$ is nonelementary.
\end{proof}

\begin{lemma}
\label{lemma:inv} 
Let $M$ be a nonelementary, orientable and hyperbolic three manifold. Then, for $n\geq 2$ the subspace of invariants of $V_n$ is trivial:
$$
V_n^{\pi_1(M)}=0.
$$
\end{lemma}

\begin{proof}
 Let us fix a basis for $V_n$. Let 
	$e_1=\left( \begin{smallmatrix} 1 \\ 0 \end{smallmatrix}\right)$ and 
	$e_2=\left( \begin{smallmatrix} 0 \\ 1 \end{smallmatrix}\right)$,
	so that $\{e_1,e_2\}$ is the standard basis for $V_2=\mathbf C^2$.
	Thus
	\[
	\{ e_1^{n-1}, e_1^{n-2} e_2,\ldots , e_2^{n-1}\}
	\]
	is a basis for $V_n = \operatorname{Sym}^{n-1}(V_2)$.

	Since $M$ is nonelementary, there exists at least one element $\gamma\in\pi_1(M)$ whose holonomy is nonparabolic   (cf. \cite[Corollary~3.25]{Kapovich}).
	Up to conjugation, it is 
	\[
	\pm \begin{pmatrix}
		\lambda & 0 \\
		0 & \lambda^{-1}
	\end{pmatrix}.
	\]
	for some $\lambda\in\mathbf C$, with $\vert\lambda\vert>1$.
	This means that the vectors $e_1$ and $e_2$ of the standard 
	basis for $\mathbf C^2$ are eigenvectors.
	Since $V_n$ is the $(n-1)$-symmetric power 
	of $\mathbf C^2$, for $n$ even the only element of $V_n$ $\gamma$-invariant is zero. 
	For $n$ odd, the subspace of $\gamma$-invariants of $ V_n$ is the line generated by $e_1^{\frac{n-1}2}e_2^{\frac{n-1}2}$. 
	Any other matrix of $\operatorname{SL}(2,\mathbf C)$ that fixes $e_1^{\frac{n-1}2}e_2^{\frac{n-1}2}$ 
	is either diagonal or antidiagonal (zero entries in the diagonal). Antidiagonal matrices
	have trace zero, hence they have order four, so they cannot occur because the holonomy of $M$ 
	has no torsion elements.
	Also, any element $\gamma'\in\pi_1(M)$ that does not commute with $\gamma$ has nondiagonal holonomy, thus $0$ is the only
	element of $ V_n$ invariant by both $\gamma$ and $\gamma'$.
\end{proof}

\begin{proof}[Proof of Proposition~\ref{prop:invariant}]

	When $F$ has genus $g\geq 2$, then by Lemma~\ref{lemma:bdry}  $\operatorname{Hol}(\pi_1(F))\backslash \mathbf H^3$ is a nonelementary hyperbolic
	3-manifold. We apply Lemma \ref{lemma:inv} to conclude that $V_n^{\pi_1(F)}=0. $
	
	Assume now that $F$ is a torus $T^2$. After conjugation, elements of $\pi_1(T^2)$ have holonomy
	\[
	\pm \begin{pmatrix} 1 & \tau \\
  		0 & 1
	\end{pmatrix} \in \operatorname{SL}(2,\mathbf C).
	\]
	The previous matrix maps $e_1^{n-i-1} e_2^i$ to $(\pm 1)^{n-1} e_1^{n-i-1} (e_2+\tau e_1)^i$,
	and it follows easily that there is no invariant subspace when $n$ is even and the lift 
	is nonpositive or it is generated by $e_1^{n-1}$ otherwise.
\end{proof}

Applying Lemma~\ref{lemma:H_F}, Proposition~\ref{prop:invariant}, 
Theorem~\ref{thm:inclusion} and Lemma~\ref{lemma:inv}, we get the following corollaries.

\begin{corollary}
	\label{cor:h1mg}
	Let $M$ be a hyperbolic manifold with $k$ cusps and 
	$l$ ends of infinite volume of genus $g_1,\dots,g_l$,
	and let $n \geq 2$. 
	Then
	\begin{eqnarray*}
		\dim_\mathbf{C} H^0( \partial \overline M; E_{\rho_n} ) & = &  a, \\
		\dim_\mathbf{C} H^1( \partial \overline M; E_{\rho_n} ) & = &  {\sum_{i=1}^l 2n(g_i-1)} + 2a, \\
		\dim_\mathbf{C} H^2( \partial \overline M; E_{\rho_n} ) & = &  a,
	\end{eqnarray*}
	where $a$ is equal to $k$ if $n$ is odd, and equals to 
	the number of cusps for which the lift of the holonomy is positive if $n$ is even. 
\end{corollary}

\begin{corollary}
	\label{cor:h1m}
Let $M$ be as in Corollary~\ref{cor:h1mg}. Then
	$H^0(   M; E_{\rho_n} )  =  0$, 
	\[
	 \dim_\mathbf{C} H^1(  M; E_{\rho_n} )  =   {\sum_{i=1}^l n(g_i-1)} + a,
	\]
	and $\dim_\mathbf{C} H^2(  M; E_{\rho_n} )  =   a$.
\end{corollary}

\subsection{Lifts of the holonomy representation}
\label{sec:lift}

\begin{proposition}[\cite{culler}]
The holonomy representation of a hyperbolic 3-manifold $M$ lifts to $\operatorname{SL}(2,\mathbf C)$.
In addition, there is a natural bijection between the set of lifts and the set of spin structures.
\end{proposition}

This is proved in Section~2 of \cite{culler}. Essentially the idea is that a spin structure on $M$ has a section, because 
$M$ is parallelizable, and this section lifts to a equivariant section of the spin bundle on the universal covering of $M$.
Identifying the universal covering of $M$ with $\mathbf H^3$, the spin bundle corresponds to $\operatorname{SL}(2,\mathbf C)$,
and equivariance of the section gives the lifted representation of $\pi_1(M)$ in 
$\operatorname{SL}(2,\mathbf C)$.  Notice that on both sets, the set of spin structures and the set of lifts, there is a simply transitive action
of $H^1(M;\mathbf Z/2\mathbf Z)$. We view elements in  $H^1(M;\mathbf Z/2\mathbf Z)$ as homomorphisms  $\pi_1(M)\to \mathbf Z/2\mathbf Z$
that
describe the difference between signs of  two different lifts.

Assume that $M$ has $k$ cusps, and choose $\gamma_1,\ldots,\gamma_k\in\pi_1(M)$ $k$ elements so that each
$\gamma_i$ is represented by a simple closed curve in one of the torus of the cusp, and different curves 
go to different cusps.

\begin{lemma} 
\label{lemma:trace-2}
For any choice of curves as above, there exists a lift 
\[
\rho\colon\pi_1(M)\to \operatorname{SL}(2,\mathbf C)
\]
of the holonomy representation such that $\operatorname{trace}(\rho(\gamma_i))=-2,$
for $i=1,\ldots,k$.
\end{lemma}

\begin{proof}
We denote the peripheral torus by $T_1^2,\ldots, T_k^2$.
Let $\mu_i\in\pi_1(T^2_i)$ be represented by a simple closed curve intersecting $\gamma_i$ in one point,
so that $\gamma_i$ and $\mu_i$ generate $\pi_1(T^2_i)$. 
We can replace $\gamma_i$ by 
$\gamma_i\mu_i^{2 n_i}$, for any integer $n_i$, 
as multiplying by an even power of $\mu_i$ does not change the sign of the trace.
We chose the $n_i$ sufficiently large so that Thurston's hyperbolic Dehn filling applies to these slopes. More precisely, we require that
 there is a continuous path of cone manifold structures with cone angle $\alpha\in [0,2\pi]$, so that $\alpha=0$ is the complete
structure on $M$ and $\alpha=2\pi$ is the filled manifold (cf. \cite{ThurstonNotes, HKAnnals}).
Now we chose the lift of the hyperbolic structure on the filled manifold,
using Culler's Theorem \cite{culler}, and consider the induced lifts  corresponding to changing continuously the cone angle.
The map $X(M,\operatorname{SL}(2,\mathbf C))\to X(M,\operatorname{PSL}(2,\mathbf C))$ is a local homeomorphism except at characters of reducible representations
or representations that preserve a (unoriented) geodesic of $\mathbf H^3$  \cite{HP2}. 
Thus we get  a continuous path of representations in $X(M,\operatorname{SL}(2,\mathbf C))$ parametrized by the cone angle
$\alpha\in [0,2\pi]$, cf. \cite[Thm.~4.1]{culler}.

The holonomy of $\gamma_i$ is conjugate to 
$$
\pm\begin{pmatrix}
    \exp(i\,\alpha/2) & 0 \\ 0 & \exp(- i\,\alpha/2)
   \end{pmatrix}
$$
and its trace is $\pm 2\cos (\alpha/2)$.  The sign $\pm$ must be constant by continuity. This is clear when
$\alpha\neq \pi$ because then the trace is nonzero. When 
$\alpha=\pi$, we use the local rigidity theorem of \cite{HKJdG,Weiss},
that says that this path is locally parametrized by $\alpha$, and since the derivative of $\pm 2\cos (\alpha/2)$
at $\alpha=\pi$ is $\pm\sin(\pi/2)=\pm 1$, the trace is monotonic on $\alpha$ when  $\alpha=\pi$.

Finally, since we have chosen a lift that is trivial on $\gamma$  when $\alpha=2\pi$, the choice of sign is 
$$
- 2\cos(\alpha/2),
$$ and when
$\alpha=0$ we get the result.
\end{proof}

We obtain the following well known result, proved by Calegari in \cite{calegari},  that applies for instance
to the longitude of a knot.

\begin{corollary}
\label{cor:tr-2}
Let $\gamma$ be a simple closed curve in a torus of $\partial\overline M$ homotopically nontrivial.
If $\gamma$ is homologous to zero in $H_1(M;\mathbf Z/2\mathbf Z)$, then, for every lift 
$\varphi\colon\pi_1(M)\to \operatorname{SL}(2,\mathbf C)$ of the holonomy representation,
$$
\operatorname{trace}(\varphi(\gamma))=-2.
$$
\end{corollary}

\begin{proof}
 The proof follows from the fact that the sign of $\varphi(\gamma)$ cannot be changed by taking different lifts, and by applying
Lemma~\ref{lemma:trace-2}.
\end{proof}

\begin{corollary}
	\label{cor:nonpositivelifts}
	Let $M$ be a hyperbolic manifold with a single cusp. Then all 
	lifts of the holonomy representation are nonpositive on $\pi_1 (\partial M)$. 
\end{corollary}

\begin{proof}
	Since the inclusion in homology 
	\[
	H_1(U; \mathbf Z/2\mathbf Z)\to H_1( M; \mathbf Z/2\mathbf Z)
	\] 
	has rank one, there exists a simple closed curve representing a nontrivial element 
	in $H_1(T^2; \mathbf Z/2\mathbf Z)\cong H_1(U; \mathbf Z/2\mathbf Z)$ that is 
	$\mathbf Z/2\mathbf Z$-homologous to zero in $M$. Thus Corollary~\ref{cor:tr-2} 
	applies here, and every lift of the holonomy restricted to the peripheral group is nonpositive. 
\end{proof}

\begin{proof}[Proof of Theorem~\ref{thm:H1cusp}] Apply Corollaries~\ref{cor:h1m} and
	\ref{cor:nonpositivelifts}.
\end{proof}

\section{Infinitesimal Rigidity}
\label{section:infinitesimalrigidity}

Here we prove Theorem~\ref{thm:l2rigidity}, that we restate.

\begin{theorem}
	\label{thm:l2rigidityagain}
	Let $M$ be a complete hyperbolic $3$-manifold that is topologically finite. 
	If $\partial \overline M$ is the union of $k$ tori and $l$ surfaces of genus $g_1,\ldots,g_l\geq 2$,
	and $n\geq 2$, then 
	\[
	\dim_\mathbf{C} H^1(M; E_{\operatorname{Ad}\circ \rho_n} ) = k(n-1) + \sum (g_i-1)(n^2-1).
	\]
	In particular, if $M$ is closed then $H^1(M; E_{\operatorname{Ad}\circ \rho_n}) = 0.$
	In addition, all nontrivial elements in $H^1(M; E_{\operatorname{Ad}\circ \rho_n} )$
	are nontrivial in $H^1(\partial\overline M,E_{\operatorname{Ad}\circ \rho_n})$ and have no $L^2$ representative.
\end{theorem}

\begin{proof} 
	By Lemma \ref{lemma:SL2Ad} we have
	$\Lie{sl}(n,\mathbf{C}) \cong V_{2n - 1} \oplus V_{2n - 3} \cdots \oplus V_3.$
	Hence, 
	\begin{equation}
		\label{eqn:h1decompose}
		H^1(M; E_{\operatorname{Ad}\circ \rho_n} ) \cong H^1(M; E_{\rho_{2 n -1}}) \oplus  
		H^1(M; E_{\rho_{2 n-3} })\oplus \cdots \oplus H^1(M; E_{\rho_3}).
	\end{equation}
	The theorem now follows from this isomorphism, Corollary \ref{cor:h1mg} and
	Theorem \ref{thm:inclusion}.
\end{proof}

Next we want to prove Theorem~\ref{thm:smoothness}. See \cite{LubotzkiMagid} for basic results
about representation and character varieties.
The variety of representations of $\pi_1(M)$ in $\operatorname{SL}(n,\mathbf C)$ is
$$
R(M,\operatorname{SL}(n,\mathbf C))=\hom(\pi_1(M), \operatorname{SL}(n,\mathbf C)).
$$
Since $\pi_1(M)$ is finitely generated, this is an algebraic affine set. 
The group $\operatorname{SL}(n,\mathbf C)$ acts by conjugation  on
$R(M,\operatorname{SL}(n,\mathbf C))$ algebraically, and the quotient in the 
algebraic category is the variety of characters:
$$
X(M,\operatorname{SL}(n,\mathbf C))=R(M,\operatorname{SL}(n,\mathbf C))/\!/ \operatorname{SL}(n,\mathbf C).
$$
For a representation $\rho\in R(M,\operatorname{SL}(n,\mathbf C)) $ its character is the map
$$
\begin{array}{rcl}
	\chi_\rho\!:\!\pi_1(M) & \to & \mathbf C\\
	\gamma & \mapsto & \operatorname{trace}(\rho(\gamma)).
\end{array}
$$
The projection $R(M,\operatorname{SL}(n,\mathbf C))\to 
X(M,\operatorname{SL}(n,\mathbf C)) $ maps each representation $\rho$ to its character $\chi_\rho$.

Weil's construction gives a natural isomorphism between the Zariski tangent space to a representation
$T^{Zar}_{\rho} R(M,\operatorname{SL}(n,\mathbf C))$ and $
Z^1(\pi_1(M),V_{\operatorname{Ad}\circ\rho})
$,
the space of group cocycles valued in
the lie algebra $\Lie{sl}(n,\mathbf C)$, which as $\pi_1(M)$-module is also written as $V_{\operatorname{Ad}\circ\rho}$.
Namely, 
$
Z^1(\pi_1(M),V_{\operatorname{Ad}\circ\rho})
$
is the set of maps $d\!:\pi_1(M)\to V_{\operatorname{Ad}\circ\rho} $ that satisfy the cocycle relation
$$
d(\gamma_1\gamma_2)= d(\gamma_1)+ \operatorname{Ad}_{\rho(\gamma_1)} d(\gamma_2),\qquad \forall \gamma_1,\gamma_2\in\pi_1(M).
$$
Notice that $ R(M,\operatorname{SL}(n,\mathbf C))$ may be a non reduced algebraic set, so the Zariski tangent space 
may be larger than the Zariski tangent space of the underlying algebraic variety.

The space of coboundaries $
B^1(\pi_1(M),V_{\operatorname{Ad}\circ\rho})
$
is the set of cocycles that satisfy $d(\gamma)= \operatorname{Ad}_{\rho(\gamma) } m- m$ for all $\gamma\in\pi_1(M)$ and for some fixed
$m\in V_{\operatorname{Ad}\circ\rho}$. 
The space of coboundaries is the tangent space to the orbit by conjugation, so
under some  hypothesis  the cohomology
may be identified with the tangent space of the variety of characters (Proposition~\ref{prop:tangentspaces}).
Since $M$ is aspherical, the group cohomology of $\pi_1(M)$
$$
H^1(\pi_1(M); V_{\operatorname{Ad}\circ\rho} )=Z^1(\pi_1(M),V_{\operatorname{Ad}\circ\rho})/B^1(\pi_1(M),V_{\operatorname{Ad}\circ\rho})
$$
is naturally isomorphic to $H^1(M; E_{\operatorname{Ad}\circ\rho})$.

\begin{definition}
A representation $\rho\!:\pi_1(M)\to \operatorname{SL}(n,\mathbf C)$ is \emph{semisimple}  if every subspace of $\mathbf C^n$
invariant by $\rho(\pi_1(M))$ has an invariant complement.
\end{definition}

Thus a semisimple representation  decomposes as direct sum of 
simple representations, where simple means without proper 
invariant subspaces.

The following summarizes the relation between tangent spaces and cohomology. See \cite{LubotzkiMagid} for a proof.

\begin{proposition} 
\label{prop:tangentspaces}
Let $\rho\in R(M,\operatorname{SL}(n,\mathbf C))$.
\begin{enumerate}
 \item There is a natural isomorphism
$$
Z^1(\pi_1(M),V_{\operatorname{Ad}\circ\rho})\cong T^{Zar}_{\rho} R(M,\operatorname{SL}(n,\mathbf C)).
$$
\item If  $\rho$ is semisimple, then it induces an isomorphism
$$
H^1(\pi_1(M); V_{\operatorname{Ad}\circ\rho})\cong T^{Zar}_{\rho} X(M,\operatorname{SL}(n,\mathbf C)).
$$
\item If  $\rho$ is semisimple and a smooth point of $R(M,\operatorname{SL}(n,\mathbf C))$, then its character
$\chi_{\rho}$ is a smooth point of $X(M,\operatorname{SL}(n,\mathbf C)$.
\end{enumerate}
\end{proposition}

A point in an algebraic affine set is smooth iff it has the same dimension that its Zariski tangent space. So 
to prove smoothness we need to compute these dimensions.

\begin{lemma}
\label{lemma:torussmooth}
Let $\rho_n$ be as in Theorem~\ref{thm:smoothness}, and $T^2$ a component of $\partial\overline M$ corresponding
to a cusp. Then the restriction of $\rho_n$ to $\pi_1(T^2)$ is a smooth point of
$R(T^2,\operatorname{SL}(n,\mathbf C))$.
\end{lemma}

\begin{proof}
Knowing  that  $\dim R(T^2,\operatorname{SL}(n,\mathbf C))\leq \dim Z^1(T^2,V_{\operatorname{Ad}\circ\rho_n})$,
we want to show that equality of dimensions holds.
Before the cocycle space, we first compute the dimension of the cohomology group. 
By Equation~\ref{eqn:h1decompose} in the proof of Theorem~\ref{thm:l2rigidityagain}:
\[
\dim H^1(T^2; E_{\operatorname{Ad}\circ\rho_n}) = \sum_{i=2}^n \dim H^1(T^2; E_{\rho_{2i-1}}).
\]
Hence, by Corollary~\ref{cor:h1mg}, 
$$
\dim H^1(T^2; E_{\operatorname{Ad}\circ\rho_n})= 2(n-1).
$$
We apply the same splitting for computing the dimension of the coboundary space. It is the sum of terms 
$\dim B^1(T^2, E_{\rho_{k}})$, for $k$ odd from $3$ to $2 n-1$.
Since we have an exact sequence
$$
0\to V_k^{\pi_1(T^2)} \to V_k\to B^1(T^2, E_{\rho_{k}})\to 0,
$$
$\dim B^1(T^2, E_{\rho_{k}}) = k - \dim V_k^{\pi_1(T^2)}= k-1$, by Lemma~\ref{prop:invariant}.
Thus 
$$
\dim B^1(T^2,E_{\operatorname{Ad}\circ\rho_n})= (2n-2)+(2n-4)+\cdots+2=n^2-n.
$$
Hence as $H^1(T^2; E_{\operatorname{Ad}\circ\rho_n})=Z^1(T^2,E_{\operatorname{Ad}\circ\rho_n})/B^1(T^2,E_{\operatorname{Ad}\circ\rho_n})$, we have:
\begin{eqnarray*}
	\dim Z^1(T^2,E_{\operatorname{Ad}\circ\rho_n}) & = & \dim H^1(T^2,E_{\operatorname{Ad}\circ\rho_n}) +
	\dim B^1(T^2,E_{\operatorname{Ad}\circ\rho_n}) \\ 
	&=& n^2+n-2.
\end{eqnarray*}

Now we look for a lower bound of $\dim R(T^2,\operatorname{SL}(n,\mathbf C))$. 
Fix $\{\gamma_1,\gamma_2\}$ a generating set of $\pi_1(T^2)$.
The representation
$\rho_n$ restricted to $\pi_1(T^2)$ has eigenvalues equal to $\pm 1$. By deforming 
the representation of $\pi_1(T^2)$ to $\operatorname{SL}(2,\mathbf C)$, and by composing it with
the representation of $\operatorname{SL}(2,\mathbf C)$ to $\operatorname{SL}(n,\mathbf C)$,  there exists a representation 
$\rho'\in R(T^2,\operatorname{SL}(n,\mathbf C)) $ arbitrarily close to $\rho_n$ such that all 
eigenvalues of $\rho'( \gamma_1)$ are different,
in particular $\rho'( \gamma_1)$ diagonalises. Now, to find deformations of $\rho'$,
notice that $\rho'( \gamma_1)$ can be deformed with $n^2-1=\dim (\operatorname{SL}(n,\mathbf C))$ 
parameters, and having all eigenvalues different is an open condition. As
$\rho'( \gamma_2)$ has to commute with $\rho'( \gamma_1)$, it has the same eigenspaces,
but one can still chose $n-1$ eigenvalues for $\rho'( \gamma_2)$. This proves that the dimension of some irreducible
component of $R(T^2,\operatorname{SL}(n,\mathbf C))$ that contains $\rho_n$ is at least
$$
n^2-1+ n-1= n^2+n-2.
$$
As this is $\dim Z^1(T^2,E_{\operatorname{Ad}\circ\rho_n})$, it is  a smooth point.
\end{proof}

\begin{proof}[Proof of Theorem~\ref{thm:smoothness}]
Using Proposition~\ref{prop:tangentspaces}, we just prove that $\rho_n$ is a smooth
point of the variety of representations.

Given a Zariski tangent vector $v\in Z^1(M, V_{\operatorname{Ad}\circ\rho_n})$, we have to show that it is integrable, 
i.e. that here is a path in the variety of representations whose tangent vector is $v$.
For this, we use the algebraic obstruction theory, see \cite{goldman, HP}. There exist an infinite sequence
of obstructions that are cohomology classes in $H^2(M, V_{\operatorname{Ad}\circ\rho_n})$, each obstruction being defined
only if the previous one vanishes. These  are related to the analytic expansion
in power series of a deformation of a representation,
and to Kodaira's theory of infinitesimal deformations. 
Our aim is to show that this infinite sequence vanishes. This gives a formal power series,
that does not need to converge, but this is sufficient for $v$ to be a tangent vector by a theorem of Artin \cite{Artin}
(see \cite{HP} for details).
We do not give the explicit construction of these obstructions, we just use that they are natural and that they live
in the second cohomology group.

By Theorem~\ref{thm:inclusion} we have an isomorphism:
\begin{equation}
 \label{eqn:isoH2}
H^2(M;E_{\operatorname{Ad}\circ\rho_n})\cong H^2(\partial\overline M;E_{\operatorname{Ad}\circ\rho_n}).
\end{equation}
Now, 
$H^2(\partial\overline M;E_{\operatorname{Ad}\circ\rho_n})$ decomposes as the sum of the connected components of $\partial \overline M$.
If $F_g$ has genus $g\geq 2$ then $H^2(F_g ;E_{\operatorname{Ad}\circ\rho_n})=0$. Thus, only the components of $\partial \overline M$ 
that are tori appear in $ H^2(\partial\overline M;E_{\operatorname{Ad}\circ\rho_n})$. 
By Lemma~\ref{lemma:torussmooth} and naturality, the obstructions vanish when restricted to $H^2(T^2;E_{\operatorname{Ad}\circ\rho_n})$,
hence they vanish in $H^2(M;E_{\operatorname{Ad}\circ\rho_n})$ by the isomorphism (\ref{eqn:isoH2}).
\end{proof}


\appendix
\section{Some results on principal bundles}
\label{principalbundles}

Throughout this section $P$ will denote a $G$-principal bundle over 
a manifold $M$.

\begin{rem}
	We will follow the convention that the action of $G$ is on the right.
\end{rem}

Assume we have a connection on $P$ with connection form $\omega \in \Omega^1(P;\Lie{g})$.
This connection defines a horizontal vector bundle $H$ on $P$. 
The differential of the bundle projection $\pi_P\colon P\rightarrow M$ is an isomorphism when restricted to $H$. 
Hence, given $X_p \in TM$ and $u\in \pi_P^{-1}(p)$, there exists a unique 
$\tilde{X}_u \in H_u$ that is projected to $X_p$. The vector $\tilde{X}_u$ is called 
the horizontal lift of $X_p$ at $u$. A vector field on $P$ is called horizontal 
if it is tangent to $H$.

All these definitions can be extended in a natural way to the cotangent bundle,
exterior powers, tensor powers, etc. Therefore, it makes sense to talk about 
horizontal forms, horizontal tensors, etc.

Let's recall a common construction. Let $F$ be a differentiable manifold on which $G$ acts on the left. 
The associated bundle, denoted by $P\times_G F$, is the quotient of $P\times V$ 
by the diagonal right action of $G$ 
(i.\;e.\;if $(u, x)\in P\times F$, then $(u, x) \cdot g = (ug,g^{-1}x)$). 
The space $P\times_G F$ has in a natural way a structure of fiber bundle over $M$ with typical
fiber $F$.

\begin{obs}
	The definition of $P\times_G F$ allows us to interpret a point
	$u$ in $P$ as an isomorphism between $F$ and the fiber of $P\times_G F$ at $\pi_P(u)$. 
	Let's say, if $\pi$ denotes the quotient map $ P\times V \rightarrow P\times_G F$, then 
	$\pi(u,\cdot )$ is an isomorphism. Note that $\pi( u\cdot g, x) = \pi(u, g x)$.
\end{obs}

We can generalize the notion of associated bundle just ``twisting $F$''; that is,
we can take as a starting point an arbitrary bundle over $P$ with typical fiber $F$, instead of just the product 
bundle $P\times F$. Let $\pi_Q\colon Q\rightarrow P$ be a bundle over $P$ with typical fiber $F$. Assume that
we have a fiber-preserving action (on the right) of $G$ on $Q$ that is compatible with the action on $P$ 
(i.\;e.\;$\pi_Q(q\cdot g) = \pi_Q(q)\cdot g$). The quotient
$Q/G$ is in a natural way a fiber bundle over $M$ with typical fiber $F$. In this case, a point $u \in P$
can be interpreted as an isomorphism between the fiber of $Q$ at $u$, and the fiber of $Q/G$ at $\pi(u)$. 

\begin{proposition}
	There is a canonical isomorphism between the space of $G$-equivariant sections of $Q$, and 
	the space of sections of the associated bundle $Q/G$.
\end{proposition}

Now we want to specialize all these things to the case $Q = \bigwedge^r H^*\otimes V$, where $V$
is a fixed vector space. Let's fix a linear representation $\rho\colon G\rightarrow \operatorname{Aut}(V)$, in such a
way that $V$ becomes a left $G$-module. We then let act $G$ on $Q$ on the right as follows: if $\alpha_p\otimes w_p$ belongs to
$Q_p$, then $(\alpha_p\otimes w_p)\cdot g = R_{g^{-1}}^*\alpha_p\otimes \rho(g)^{-1}w_p \in Q_{pg}$. 
Using horizontal lifts we can identify $Q/G$ with $\bigwedge^r T^*M \otimes E$. More precisely, let $p\in M$, $u\in \pi^{-1}(p)$, and 
$\operatorname{H}_u\colon T_p M\rightarrow H_u$ the horizontal lift map. Then, if we interpret $u$ an isomorphism 
between $V$ and $E_p$, we obtain the isomorphism 
$\varphi_u\colon\operatorname{H}^*_u\otimes u\colon Q_u\rightarrow \bigwedge^r T_p*M \otimes E_p $.
Since horizontal lift and $u$ commute with the action of $G$, we have $\varphi_{u}(v) =\varphi_{ug}(vg)$, for 
all $v\in Q$. Therefore, we get an isomorphism $\varphi$ between $Q/G$ and $\bigwedge^r T^* M \otimes E$.
 
We will denote by $\Omega_{\operatorname{Hor}}^*(P; V)^G$ the space of horizontal $V$-valued differential 
forms over $P$ that are $G$-equivariant, or, equivalently, the space of $G$-equivariant sections of the bundle 
$\bigwedge^r H^*\otimes V$. 

\begin{obs}
	A form $\alpha$ is horizontal if, and only if, it vanishes on vertical directions, 
	that is, $i_X \alpha = 0$ for any vertical vector field $X$. Also,
	$\alpha$ is $G$-equivariant if, and only if, $R_g^* \alpha = \rho(g^{-1}) \alpha$ for all $g\in G$. Therefore,
	$\alpha\in \Omega^r (P; V)$ belongs to $\Omega_{\operatorname{Hor}}^r(P;V)^G$ if, and only if, 
	\begin{eqnarray}
		R_g^*\alpha&=&\rho(g)^{-1}\alpha,\quad \text{per a tot $g\in G$},\\
		i_{Y}\alpha&=&0,\quad \text{per a tot $Y\in\Lie{g}$.}
	\end{eqnarray}
	Note that we are identifying $\Lie{g}$ with the space of $G$-invariant vertical vectors over a fixed fiber of $P$.
\end{obs}

The connection on $P$ defines an exterior covariant differential on $G$-equivariant horizontal forms. Namely,
$$ D\alpha = (d \alpha)\circ \pi_h,\quad \mbox{for} \quad \alpha\in \Omega^r_{\operatorname{Hor}}(P; V)^G$$
where $\pi_h$ is projection on the horizontal distribution defined by the connection. On the other hand, a connection 
on $P$ induces a connection on the vector bundle $P\times_\rho V$, and hence an exterior covariant differential $d_\rho$ on
$\Omega^r(M; E)$. It is easily verified that the canonical isomorphism between the spaces $\Omega^*_{\operatorname{Hor}}(P; V)^G$ and
$\Omega^*(M; E)$, ``commute'' with exterior covariant differentiation (see \cite[p. 76]{Kob}).

\begin{proposition}
	\label{prop:dif_cov}
	Let $\omega\in \Omega^1(P;\Lie{g})$ be the connection form of the connection 
	defined on $P$. Then the following formula holds 
	$$D\alpha = d\alpha + \rho(\omega)\wedge \alpha.$$
\end{proposition}

\begin{rem}
	If $V_1,\dots,V_{p+1}$ are vector fields on $P$, by definition,
	$$(\rho\wedge\alpha)(V_1,\dots,V_{p+1}) = 
	\sum_{i=1}^{p+1} (-1)^{i+1}\rho(Y_i)\big(\alpha(V_1,\dots,\widehat{V_i},\dots,V_{p+1})\big).$$
	Taking a base of $V$, $\rho(w)$ is just a matrix of $1$-forms, $\alpha$ a column vector of 
	$p$-forms, and the product $\rho(\omega)\wedge \alpha$ is just the product of a matrix by 
	a vector.

\end{rem}

\begin{proof}
	We must prove the form $d\alpha + \rho(\omega)\wedge \alpha$ is horizontal, and that on horizontals vectors
	coincides with $D\alpha$. The second fact is obvious from the the definition of $D$ and the fact that 
	$\omega$ vanishes on horizontal vectors. Hence we only need to prove that $d\alpha + \rho(\omega)\wedge \alpha$
	vanishes on vertical vectors. Let be $X^*$ the fundamental vector field associated to $X\in \Lie{g}$, using Cartan's identity
	($L_X^* = d i_X^* + i_X^* d)$) we get
	$i_X^*(d\alpha + \rho(\omega)\wedge \alpha) = L_X^* \alpha - d(i_X^*\alpha) + \rho(X)\alpha).$
	The infinitesimal version of the $G$-equivariance of $\alpha$ states that $L_X^* \alpha = -\rho(X)\alpha$. Then we conclude
	that $d\alpha + \rho(\omega)\wedge \alpha$ is vertical.
\end{proof}

Now assume that $M$ is a Riemannian manifold, and that we have a metric on the vector 
bundle $E = P\times_G V$. These metrics induce an inner product on the space of 
$E$-valued forms over $M$.
$$(\alpha, \beta) = \int_M \langle\alpha(x),\beta(x)\rangle_x \omega_M.$$

On the other hand, the Riemannian metric on $M$ defines a metric tensor on the 
horizontal bundle $H$, in such a way that horizontal lifts are isometries. Also, the metric defined on $E$
defines a metric on the trivial vector bundle $P\times V$. A right invariant volume form $\omega_G$ on $G$ 
defines a right invariant volume form along the fibers of $P$. Therefore, we can define an inner 
product on $\Omega_{\operatorname{Hor}}^r(P; V)^G$ by
$$(\tilde{\alpha}, \tilde{\beta}) = \int_P \langle \tilde{\alpha}(u), \tilde{\beta}(u) \rangle_u \pi_P^*(\omega_M)\wedge\omega^*_G.$$

We want to study how the metrics defined on $\Omega^r(M; E)$ and $\Omega_{\operatorname{Hor}}^r(P; V)^G$ are related by
the canonical isomorphism. However, this comparison doesn't make sense if $G$ is not assumed to be compact
(if $\alpha\in\Omega^r(M; E)$ has compact support, then the corresponding form $\tilde{\alpha}$ 
in $\Omega_{\operatorname{Hor}}^r(P; V)^G$ has compact support if, and only if, $G$ is compact). 
From now on we will assume that $G$ is compact. In order to avoid confusions we will denote $G$ by $K$ in this case.
In this case we can simplify things a little bit. First, take a $K$-invariant metric on $V$, and use it to define a ``constant'' metric on $P\times V$. Since this metric is $K$-invariant, we get a metric
on the vector bundle $E$. Under these hypothesis, we get a nice relation between these two metrics. In order to get this relation, we need the following lemma.

\begin{proposition}
	Let $\omega_M$ be a volume form on $M$, and $\omega_K$ a right invariant volume form on $K$.
	Denote by $\omega^*_K$ the right invariant volume form on the fibers of $P$ defined by $\omega_K$.
	If $f$ is a function defined on $P$, then the function 
	$\bar{f}(u) = \int_K f(ug) \omega_K$ is invariant along the fibers, and hence can be seen as a function on $M$.
	With these hypothesis, we have
	$$ \int_P f(u) \pi_P^*(\omega_M)\wedge\omega_K^* = \int_M \bar{f}(x) \omega_M$$    
\end{proposition}
\begin{proof}
	Take an open set $U\subset M$ that trivializes $P$, and a trivializing map 
	$\psi\colon U\times K\rightarrow \pi_P^{-1}(U)$. Let's denote by $\pi_U$ and
	$\pi_K$ the projection of $U\times K$ on the first and on the second factor
	respectively. We have $\omega^*_K = (\psi^*)^{-1}(\pi_K^*(\omega_K))$. The change
	of variable formula gives
	$$\int_{\pi_P^{-1}(U)} f(u) \pi_P^*(\omega_M)\wedge\omega_K^* = 
	\int_{U\times K} f(\psi(x, g)) \pi_U^*(\omega_M)\wedge\pi_K^*(\omega_K).$$
	By Fubini's Theorem, the last integral is,
	$$\int_U \left( \int_K f(\psi(x, g)) \omega_K \right)\omega_M = \int_U \bar{f}(x) \omega_M.$$
	The result follows by taking a partition of unity subordinated to a trivializing open cover.
\end{proof}

The function $\langle \tilde{\alpha}(u), \tilde{\beta}(u) \rangle_V$ is constant along the fibers, and equals to
$\langle\alpha(x),\beta(x)\rangle_x$, where $x = \pi_P(u)$. The above lemma then implies the following proposition.
\begin{proposition}
	\label{G_bundle:int}
	With the above notation, $$(\tilde{\alpha}, \tilde{\beta}) = \mu(K) (\alpha, \beta),$$ 
	where $\mu$ denotes the measure defined by the volume form $\omega_K$.
\end{proposition}

Consider the pairing
$$
\begin{array}{ccc}
	\Omega_{\operatorname{Hor}}^r(P;V)^K \times  \Omega_{\operatorname{Hor}}^{m-r}(P;V)^K & \longrightarrow & \mathbf{R}\\
	(\alpha, \beta) &\longmapsto&\int_P (\alpha\wedge\beta)\wedge\omega_K, 
\end{array}
$$
where the wedge product of a $V$-valued is defined using the usual wedge product on scalar-valued forms, and the inner product 
on $V$. On the other hand, the metric on the horizontal bundle, and the orientation 
we have on it, allow us to define a Hodge star operator on the space of horizontal forms,
$$*\colon \Omega_{\operatorname{Hor}}^r(P;V)^K\longrightarrow \Omega_{\operatorname{Hor}}^{m-r}(P;V)^K.$$
Note that we have $(\alpha,\beta) =\phi(\alpha, *\beta) $

\begin{proposition}
	\label{G_bundle:adjoint}
	Let $T\colon\Omega_{\operatorname{Hor}}^r(P;V)^K\rightarrow \Omega_{\operatorname{Hor}}^{r+k}(P;V)^K$ be a
	linear operator that decreases supports. Assume we have a linear operator
	$$
	S\colon \Omega_{\operatorname{Hor}}^{m-(r+k)}(P;V)^K\rightarrow \Omega_{\operatorname{Hor}}^{m-r}(P;V)^K
	$$
	such that $\phi(T\alpha,\beta) = \phi(\alpha, S\beta)$.
	Then, the formal adjoint of $T$ is 
	$$
	T^* = (-1)^{r(m-r)}*S*\colon\Omega_{\operatorname{Hor}}^{r+k}(P;V)^K\rightarrow \Omega_{\operatorname{Hor}}^r(P;V)^K.
	$$
\end{proposition}
\begin{proof}
	Let's denote $\Omega_{\operatorname{Hor}}^r(P;V)^K$ by $M_r$.
	We have the following commutative diagram,
	$$
	\xymatrix{ 
		M_{r+k}^*\ar[r]^{T^t} & M_r^*\\
		M_{r+k}\ar[u] \ar[r]^{T^*} & M_r\ar[u],\\ 
		}
	$$
	where the vertical arrows are the isomorphisms given by the metrics, $T^t$ is the dual map of $T$, and $T^*$
	its adjoint. We can factor the metric isomorphism as $\phi(*,)$. We have the following commutative diagram
	$$
	\xymatrix{ 
		M_{r+k}^*\ar[r]^{T^t} & M_r^*\\
		M_{m-(r+k)}\ar[u]^{\phi(\cdot, )} \ar[r]^S & M_{m -r}\ar[u]^{\phi(\cdot, )}\\
		M_{r+k}\ar[u]^{*}\ar[r]^{T^*} & M_r \ar[u]^{*} 
	}.
	$$
	The proposition now follows from the fact that on degree $r$ we have $*^{-1} = (-1)^{r(m-r)}*$.

\end{proof}

\nocite{*}
\bibliographystyle{abbrv}

\end{document}